\documentclass[11pt,british]{amsart}

\newif\ifpdf
\ifx\pdfoutput\undefined\pdffalse\else\pdfoutput=1\pdftrue\fi\newif\pdf
\ifpdf\relax\else
\usepackage[T1]{fontenc}
\fi

\usepackage{babel,eucal,url,amssymb,stmaryrd,booktabs,%
enumerate,amscd,paralist}

\usepackage[pagebackref]{hyperref}

\theoremstyle{plain}
\newtheorem{thm}{Theorem}[section]
\newtheorem{lem}[thm]{Lemma}
\newtheorem{pro}[thm]{Proposition}
\newtheorem{co}[thm]{Corollary}
\newtheorem*{main}{Main Theorem}

\theoremstyle{definition}
\newtheorem{defn}[thm]{Definition}

\theoremstyle{remark}
\newtheorem{rem}[thm]{Remark}

\newtheorem*{ack}{Acknowledgements}

\newcommand{\Lie}[1]{\operatorname{\textsl{#1}}}
\newcommand{\lie}[1]{\operatorname{\mathfrak{#1}}}
\newcommand{\GL}{\Lie{GL}}

\newcommand{\SO}{\Lie{SO}}

\newcommand{\Gtwo}{\ifmmode{{\rm G}_2}\else{${\rm G}_2$}\fi}

\newcommand{\LC}{{\nabla^g}}

\newcommand{\Nt}{\tilde\nabla}
\newcommand{\Rt}{\tilde R}

\DeclareMathOperator{\tr}{tr}

\newcommand{\Hodge}{\mathord{\mkern1mu *}}
\newcommand{\hook}{\mathbin{\lrcorner}}

\newcommand{\norm}[1]{\left\lVert #1\right\rVert}

\date{\today}
\begin{document}

\title[On the geometry of closed $\Gtwo$-structures]%
{On the geometry of closed $\mathbf{G_2}$-structures}

\author{Richard Cleyton}
\address[Cleyton]{Department of Mathematics, Surge 222\\
University of California Riverside\\ CA 92521 - USA}
\email{cleyton@math.ucr.edu}

\author{Stefan Ivanov}
\address[Ivanov]{University of Sofia "St. Kl. Ohridski"\\
Faculty of Mathematics and Informatics,\\
Blvd. James Bourchier 5,\\
1164 Sofia, Bulgaria}
\email{ivanovsp@fmi.uni-sofia.bg}

\begin{abstract}
  We give an answer to a question posed recently by R.Bryant, namely
  we show that a compact $7$-dimensional manifold equipped with a
  $G_2$-structure with closed fundamental form is Einstein if and only
  if the Riemannian holonomy of the induced metric is contained in
  $G_2$. This could be considered to be a $G_2$ analogue of the
  Goldberg conjecture in almost K{\"a}hler geometry. The result was
  generalized by R.L.Bryant to closed $G_2$-structures with too
  tightly pinched Ricci tensor.  We extend it in another direction
  proving that a compact $G_2$-manifold with closed fundamental form
  and divergence-free Weyl tensor is a $G_2$-manifold with parallel
  fundamental form.  We introduce a second symmetric Ricci-type tensor
  and show that Einstein conditions applied to the two Ricci tensors
  on a closed $G_2$-structure again imply that the induced metric has
  holonomy group contained in $G_2$.
\end{abstract}

\maketitle
\setcounter{tocdepth}{2}
\tableofcontents

\section{Introduction}

A $7$-dimensional Riemannian manifold is called a $G_2$-manifold if its
structure group reduces to the exceptional Lie group $G_2$.  The existence
of a $G_2$-structure is equivalent to the existence of a non-degenerate
three-form on the manifold, sometimes called the fundamental form of the
$G_2$-manifold.  From the purely topological point of view, a
$7$-dimensional paracompact manifold is a $G_2$-manifold if and only if it
is an oriented spin manifold~\cite{LM}.


In~\cite{FG}, Fernandez and Gray divide $G_2$-manifolds into 16
classes according to how the covariant derivative of the
fundamental three-form behaves with respect to its decomposition
into $G_2$ irreducible components (see also~\cite{CS}).  If the
fundamental form is parallel with respect to the Levi-Civita
connection then the Riemannian holonomy group is contained in
$G_2$, we will say that the $G_2$-manifold or the $G_2$-structure
on the manifold is \emph{parallel}. In this case the induced
metric on the $G_2$-manifold is Ricci-flat, a fact first observed
by Bonan~\cite{Bo}.  It was shown by Gray~\cite{Gr} (see
also~\cite{Br,Sal}) that a $G_2$-manifold is parallel precisely
when the fundamental form is harmonic.  The first examples of
complete parallel $G_2$-manifolds were constructed by Bryant and
Salamon~\cite{BS,Gib2}.  Compact examples of parallel
$G_2$-manifolds were obtained first by Joyce~\cite{J1,J2,J3} and
recently by Kovalev~\cite{Kov}. Compact parallel $G_2$-manifold
will be refered to as \emph{Joyce spaces}. Examples of
$G_2$-manifolds in other Fernandez-Gray classes may be found
in~\cite{Fer,CMS}.  A central point in our argument is that the
Riemannian scalar curvature of a $G_2$-manifold may be expressed
in terms of the fundamental form and its derivatives and
furthermore the scalar curvature carries a definite sign for
certain classes of $G_2$-manifolds~\cite{FI1,Br1}.

The geometry of $G_2$-structures has also attracted much attention
from physicists.  The central issue in physics is that connections
with holonomy contained in $G_2$ plays a r{\^o}le in string
theory~\cite{GKMW,Gib1,Le1,GMW,II}. The $G_2$-connections
admitting three-form torsion have been of particular interest.

In the present paper we are interested in the geometry of closed
$G_2$-structures i.e., $G_2$-manifolds with closed fundamental
form (sometimes these spaces are called almost $G_2$-manifolds or
calibrated $G_2$-manifolds). In the sense of the Fernandez-Gray
classes, this is complementary to the physicists requirement of
three-form torsion~\cite{FI}.  Compact examples of closed
$G_2$-manifolds were presented by Fernandez~\cite{Fer}.
Supersymmetric string solutions on closed $G_2$-manifolds were
investigated in \cite{Gib1}, topological quantum field theory on
closed $G_2$-manifolds were discussed in \cite{Le1}. Robert Bryant
shows in~\cite{Br1} that if the scalar curvature of a closed
$G_2$-structure is non-negative then the $G_2$-manifold is
parallel.  The question whether there are closed $G_2$-structure
which are Einstein but not Ricci-flat then naturally arises. We
investigate this question of Bryant in the compact and in the
non-compact cases.

In the first version of the present article \cite{CI} we answered
negatively to the Bryant's question, namely, we proved that there
are no closed Einstein $G_2$- structures (other than the parallel
ones) on compact 7-manifold. In \cite{Br2} R.L.Bryant generalized
this non-existence result for closed $G_2$-structures on compact
7-manifold whose Ricci tensor is too tightly pinched.

In the present article we obtain non-existence result involving third
derivatives of the fundamental form. Namely, we prove the following
\begin{main}
A compact $G_2$-manifold with closed fundamental form and harmonic
Weyl tensor (divergence-free Weyl tensor) is a Joyce space.
\end{main}
The second Bianchi identity leads to
\begin{co}
A compact $G_2$-manifold with closed fundamental form and harmonic
curvature (divergence-free curvature tensor) is a Joyce space.
\end{co}
\begin{co}
A compact Einstein $G_2$-manifold with closed fundamental form is
a Joyce  space.
\end{co}
The latter may be considered to be a $G_2$ analogue of the
Goldberg conjecture in almost K\"ahler geometry (see e.g.
\cite{Gib1}).

The representation theory of $G_2$ gives rise to a second
symmetric Ricci type tensor on $G_2$-manifolds. Therefore one may
consider two complementary Einstein equations.  We find a
connection between the two Ricci tensors and show in
Theorem~\ref{th1}, with no compactness assumption, that if both
Einstein conditions hold simultaneously on a $G_2$-manifold with
closed fundamental form then the fundamental form is parallel.

Our main tool is the canonical connection of a $G_2$-structure and its
curvature. We will show that the Ricci tensor of the canonical connection
is proportional to the Riemannian Ricci tensor.  This leads to the
corollary that a compact $G_2$-manifold with closed fundamental form which
is Einstein with respect to the canonical connection is a Joyce space.

Our main technical tool is an integral formula which holds on any
compact $G_2$-manifold with closed fundamental form. We derive the
Main Theorem as a consequence of a more general result
,Theorem~\ref{gen}, which shows that the vanishing of the
$\Lambda^2_7$-part of the divergence of the Weyl tensor implies that a
closed $G_2$-structure is parallel on a compact 7-manifold.
\begin{ack}
The authors wish to thank Andrew Swann for useful discussions and
remarks.  This research was supported by the Danish Natural
Science Research Council, Grant 51-00-0306.  The authors are
members of the EDGE, Research Training Network HPRN-CT-2000-00101,
supported by the European  Human Potential Programme.  The final
part of this paper was done during the visit of S.I.  at the Abdus
Salam International Centre for Theoretical Physics, Trieste,
Italy.  S.I. thanks the Abdus Salam ICTP for providing support
and an excellent research  environment.  The research of S.I. is
partially supported by Contract MM 809/1998 with the Ministry of
Science and Education of Bulgaria, Contract 586/2002 with the
University of Sofia "St. Kl. Ohridski".
\end{ack}

\section{General properties of $G_2$-structures}

We recall some notions of $G_2$ geometry.  Endow ${\mathbb R}^7$ with its
standard orientation and inner product.  Let $e_1,\dots,e_7$ be an oriented
orthonormal basis which we identify with the dual basis via the inner
product.  Write $e_{i_1 i_2\dots i_p}$ for the monomial $e_{i_1} \wedge
e_{i_2} \wedge \dots \wedge e_{i_p}$.  We shall omit the $\sum$-sign
understanding summation on any pair of equal indices.

Consider the three-form $\omega$ on ${\mathbb R}^7$ given by
\begin{equation}
  \omega =e_{124} + e_{235} + e_{346}+e_{457} + e_{561} + e_{672} +
  e_{713}.\label{11}
\end{equation}
The subgroup of $\GL(7)$ fixing $\omega$ is the exceptional Lie group
$G_2$.  It is a compact, connected, simply-connected, simple Lie subgroup
of $\SO(7)$ of dimension 14~\cite{Br}.

The Hodge star operator supplies the 4-form $\Hodge\omega$ given by
\begin{equation}
  \Hodge\omega = - e_{3567} - e_{4671} - e_{5712}-e_{6123} - e_{7234} -
  e_{1345} - e_{2456}.\label{12}
\end{equation}
We let the expressions
\begin{gather*}
  \omega = \frac{1}{6}\omega_{ijk}e_{ijk},\\
  \Hodge\omega = \frac{1}{24}\omega_{ijkl}e_{ijkl}
\end{gather*}
define the symbols $\omega_{ijk}$ and $\omega_{ijkl}$. We then obtain the
following set of formulae
\begin{gather}
  \omega_{ipq}\omega_{jpq} = 6\delta_{ij}, \nonumber\\
  \omega_{ipq}\omega_{jkpq} = -4\omega_{ijk},\nonumber\\
  \omega_{ijp}\omega_{klp} = -\omega_{ijkl} + \delta_{ik}\delta_{jl} -
  \delta_{il}\delta_{jk},\label{0} \\
  \omega_{ijpq}\omega_{klpq} = -2\omega_{ijkl} + 4(\delta_{ik}\delta_{jl} -
  \delta_{il}\delta_{jk}),\nonumber\\
  \omega_{ijp}\omega_{klmp} = \delta_{ik}\omega_{jlm} -
  \delta_{jk}\omega_{ilm} + \delta_{il}\omega_{jmk} -
  \delta_{jl}\omega_{imk} + \delta_{im}\omega_{jkl} -
  \delta_{jm}\omega_{ikl}\nonumber.
\end{gather}

\begin{defn}
  A {\it $G_2$-structure} on a $7$-manifold $M$ is a reduction of the
  structure group of the tangent bundle to the exceptional group $G_2$.
  Equivalently, there exists a nowhere vanishing differential three-form
  $\omega$ on $M$ and local frames of the cotangent bundle with respect to
  which $\omega$ takes the form~\eqref{11}.  The three-form $\omega$ is
  called the {\it fundamental form} of the $G_2$-manifold $M$~\cite{Bo}.

  We will say that the pair $(M,\omega)$ is a \emph{$G_2$-manifold} with
  \emph{$G_2$-structure} (determined by) \emph{$\omega$}.
\end{defn}
\begin{rem}
  Alternatively, a $G_2$-structure can be described by the existence of a
  two-fold vector cross product $P$ on the tangent spaces of $M$.

  The fundamental form of a $G_2$-manifold determines a metric through
  $g_{ij}=\frac16\omega_{ikl}\omega_{jkl}$. This is refered to as the
  metric induced by $\omega$. We write $\LC$ for the associated Levi-Civita
  connection, $||.||^2$ for the tensor norm with respect to $g$. In addition we
  will freely identify vectors and co-vectors via the induced metric $g$.
\end{rem}

Let $(M,\omega)$ be a $G_2$-manifold.  The action of $G_2$ on the tangent
space induces an action of $G_2$ on $\Lambda^k(M)$ splitting the exterior
algebra into orthogonal subspaces, where $\Lambda^k_l$ corresponds to an
$l$-dimensional $G_2$-irreducible subspace of $\Lambda^k$:
\begin{equation*}
  \Lambda^1(M)=\Lambda^1_7,
  \quad \Lambda^2(M) = \Lambda^2_7\oplus \Lambda^2_{14}, \quad
  \Lambda^3(M)=\Lambda^3_1\oplus\Lambda^3_{7}\oplus\Lambda^3_{27},
\end{equation*}
where
\begin{gather*}
  \Lambda^2_7 = \{\alpha \in \Lambda^2(M) \vert
  \Hodge(\alpha\wedge\omega)=- 2\alpha\}, \\
  \Lambda^2_{14} = \{\alpha\in \Lambda^2(M) \vert
  \Hodge(\alpha\wedge\omega) = \alpha\}\\
  \Lambda^3_1=\{t.\omega\vert~t\in R\}, \\
  \Lambda^3_7 = \{\Hodge(\beta\wedge\omega) \vert ~\beta \in
  \Lambda^1(M)\},\\
  \Lambda^3_{27} = \{\gamma \in \Lambda^3(M) \vert ~\gamma\wedge\omega=0,
  ~\gamma\wedge\Hodge\omega = 0\}.
\end{gather*}
The Hodge star $\Hodge$ gives an isometry between $\Lambda^k_l$ and
$\Lambda^{7-k}_l$.

More generally, $V^d_{(\lambda_1,\lambda_2)}$ will denote the $G_2$
representation of highest weight $(\lambda_1,\lambda_2)$ of dimension $d$.
Note that $V^1_{(0,0)} \cong \Lambda^3_1 \cong \Lambda^4_1$ is the trivial
representation, $\Lambda^1_7 \cong V^7_{(1,0)}$ is the standard
representation of $G_2$ on $\mathbb R^7$, and the adjoint representation is
$\lie g_2 \cong V^{14}_{(0,1)} \cong \Lambda^2_{14}$.  Also note that
$V^{27}_{(2,0)} \cong \Lambda^3_{27} \cong \Lambda^4_{27}$ is isomorphic to
the space of traceless symmetric tensors $S^2_0 V^7$ on $V^7_{(1,0)}$.

\section{Ricci tensors on $G_2$-manifold}

Let $(M,\omega)$ be a $G_2$-manifold with fundamental form
$\omega$. Let $g$ be the associated Riemannian metric.
\begin{equation*}
  R_{X,Y}=[\LC_X,\LC_Y]-\LC_{[X,Y]}
\end{equation*}
is then the curvature tensor of the Levi-Civita connection $\LC$ of the
metric $g$.  The Ricci tensor $\rho$ is defined as usual as the contraction
$\rho_{ij}=R_{sijs}$, where $R_{sijs}$ are the components
\begin{equation*}
  R_{sijk}:=g(R(e_s,e_i)e_j,e_k)
\end{equation*}
of the curvature tensor with respect to an orthonormal basis
$e_1,\dots,e_7$.

\begin{defn}
On $(M,\omega)$ we may define a second symmetric tensor $\rho^\star$ by
\begin{equation}\label{1}
  \rho^\star_{sm}:=R_{ijkl}\omega_{ijs}\omega_{klm}.
\end{equation}
We will call the $\rho^\star$ the {\it $\star$-Ricci tensor} of the
$G_2$-manifold.
\end{defn}

The two Ricci tensors have common trace in the following sense.  Let
$s=\tr_g\rho=\rho_{ii}$ be the scalar curvature and let the trace of
$\rho^\star$ be denoted by $s^\star=\tr_g\rho^\star=\rho^\star_{ii}$.
\begin{pro}\label{pr}
On a $G_2$-manifold we have $s^\star=-2s$.
\end{pro}
\begin{proof}
  Apply~\eqref{0} to the definition of $s^\star$ and use
  skew-symmetry of $\Hodge\omega$ and the Bianchi identity to conclude
  that $R_{ijkl}\omega_{ijkl}=0$.
\end{proof}
\begin{defn}
  We shall use the term {\it $\star$-Einstein} for $G_2$-manifold
  $(M,\omega)$ when the traceless part of the $\star$-Ricci tensor vanishes,
  i.e., when the equation
  \begin{equation*}
    \rho^\star=\frac{s^\star}{7}g
  \end{equation*}
  holds.
\end{defn}
We define associated Ricci three-forms by
\begin{gather*}
  \begin{split}
    \rho^\star_{ijk}:=R_{ijlm}\omega_{lmk}+R_{jklm}\omega_{lmi}
    +R_{kilm}\omega_{lmj}, \\
    \rho_{ijk}:=\rho_{is}\omega_{sjk}+ \rho_{js}\omega_{ski} +
    \rho_{ks}\omega_{sij}.
  \end{split}
\end{gather*}
In terms of the Ricci forms, the Weitzenb{\"o}ck formula for the
fundamental form can be written as follows
\begin{pro}
On any $G_2$-manifold the following formula holds
\begin{equation}
  \label{w}
  d\delta\omega + \delta d\omega = \LC^*\LC\omega + \rho +
  \rho^\star.
\end{equation}\qed
\end{pro}

\section{Closed $G_2$-structures}

Let $(M^7,\omega)$ be a $G_2$-manifold with closed fundamental
form.  The two-form $\delta\omega$ then takes values in
$\Lambda^2_{14}$ \cite{Br1}.  As a consequence we get
\begin{pro}
  The following formulas are valid on a closed $G_2$-structure:
\begin{equation}
  \label{bian}
  \delta\omega_{ij}\omega_{ijk} = 0, \qquad
  \delta\omega_{ip}\omega_{pjk} + \delta\omega_{jp}\omega_{pki} +
  \delta\omega_{kp}\omega_{pij} = 0.
\end{equation}\qed
\end{pro}
It is well-known~\cite{Gr} that a $G_2$-structure is parallel if and only
if it is closed and co-closed, $d\omega=\delta\omega=0$.  The two-form
$\delta\omega$ thus may be interpreted as the deviation of $\omega$ from a
parallel $G_2$-structure.  We are going to find explicit formulae for the
covariant derivatives of the fundamental form of a closed $G_2$-structure
in terms of $\delta\omega$ and its derivatives.

\begin{defn}
  The \emph{canonical connection} $\Nt$ of a closed $G_2$-structure may be
  defined by the equation
  \begin{equation}
    \label{nat}
    g(\Nt_XY,Z)=g(\LC_XY,Z) - \frac16\delta\omega(X,e_i)\omega(e_i,Y,Z)
  \end{equation}
  for vector fields $X,Y,Z$.
\end{defn}

Using~\eqref{bian} it is easy to see that $\Nt$ is a metric
$G_2$-connection, i.e., it satisfies
\begin{equation*}
  \Nt\omega =0,\qquad \Nt g = 0.
\end{equation*}
The torsion $T$ of $\Nt$ is determined by
\begin{equation*}
  g(T(X,Y),Z) = \frac{1}{6}\delta\omega(Z,e_i)\omega(e_i,X,Y).
\end{equation*}
On a compact $G_2$-manifold the canonical connection may be characterized
as the unique $G_2$ connection of minimal torsion with respect to the
$L^2$-norm on $M$.  It may also be described by the fact that the difference
$\LC-\Nt$ takes values in $\Lambda^2_7$, the orthogonal complement of $\lie
g_2\subset\Lambda^2$ with respect to the metric induced by $g$.

From the properties of the canonical connection and $\delta\omega$ one
derives
\begin{pro}
  For a closed $G_2$-structure the following relations hold:
  \begin{gather}
    \label{der1}
    \LC_i\omega_{jkl} = \frac{1}{2}\delta\omega_{ip}\omega_{pjkl},\\
    \label{der2}
    \LC_i\omega_{jklm} =
    -\frac{1}{2}\left(\delta\omega_{ij}\omega_{klm} -
      \delta\omega_{ik}\omega_{lmj} + \delta\omega_{il}\omega_{mjk} -
      \delta\omega_{im}\omega_{jkl}\right)
  \end{gather}
  and
  \begin{equation}\label{lapl}
    \LC^*\LC\omega_{jkl} =
    \frac{1}{4}\norm{\delta\omega}^2\omega_{jkl} -
    \frac{1}{4}\left(\delta\omega_{ip}\delta\omega_{ij}\omega_{pkl} +
      \delta\omega_{ip}\delta\omega_{ik}\omega_{plj} +
      \delta\omega_{ip}\delta\omega_{il}\omega_{pjk}\right).
  \end{equation}
  \qed
\end{pro}

Applying~\eqref{nat} and~\eqref{der1} we get that the curvature
$\Rt$ of the canonical connection $\Nt$ is
related to the curvature of the Levi-Civita connection by:
\begin{gather}
  \label{curv}
  R_{ijkl} = \Rt_{ijkl} +
  \frac{1}{6}\left[\LC_i\delta\omega_{jp} -
  \LC_j\delta\omega_{ip}\right]\omega_{pkl} +
  \frac{1}{9}\delta\omega_{is}\delta\omega_{jp}\omega_{spkl} -
  \frac{1}{36}\left[\delta\omega_{ik}\delta\omega_{jl} -
  \delta\omega_{il}\delta\omega_{jk} \right]
\end{gather}

\section{Curvature of closed ${G_2}$-structures}

From here on $(M^7,\omega)$ will be a $G_2$-manifold with closed
$G_2$-structure.  We have
\begin{pro}\label{pro1}
  The Ricci tensors of a closed $G_2$-structure $(M,\omega)$ are given by
  \begin{gather}
    \label{ric1}
    \rho_{lm} = - \frac{1}{4}d\delta\omega_{sjm}\omega_{sjl} +
    \frac{1}{2}\delta\omega_{lj}\delta\omega_{mj};\\
    \label{ric2}
    \rho^\star_{lm}=d\delta\omega_{sjm}\omega_{sjl} +
    \delta\omega_{lj}\delta\omega_{mj} -
    \frac{1}{2}\norm{\delta\omega}^2\delta_{ml}.
  \end{gather}
\end{pro}
\begin{proof}
  The Ricci identities for $\omega,\Hodge\omega$ together with~\eqref{der1}
  and~\eqref{der2} lead to the following useful
 \begin{lem}
   If $\omega$ is a closed $G_2$-structure on $M^7$ then
   \begin{multline}\label{ricid1}
     \rho_{sr}\omega_{rkl} + \frac{1}{2}R_{skir}\omega_{lir} -
     \frac{1}{2}R_{slir}\omega_{kir} = \\
     - \frac{1}{4}\left(d\delta\omega_{sjp} +
       \LC_s\delta\omega_{jp}\right)\omega_{pjkl} +
     \frac{1}{2}\delta\omega_{pj}\delta\omega_{sj}\omega_{klp}.
   \end{multline}
   \begin{multline}\label{ricid2}
     -R_{sijr}\omega_{rklm} - R_{sikr}\omega_{jrlm} -
     R_{silr}\omega_{jkrm} - R_{simr}\omega_{jklr} = \\
     \frac{1}{2}\left[\left(\LC_i\delta\omega_{sj} -
         \LC_s\delta\omega_{ij}\right)\omega_{klm} -
       \left(\LC_i\delta\omega_{sk} -
         \LC_s\delta\omega_{ik}\right)\omega_{lmj}\right]\\
     + \frac{1}{2}\left[\left(\LC_i\delta\omega_{sl} -
         \LC_s\delta\omega_{il}\right)\omega_{mjk} -
       \left(\LC_i\delta\omega_{sm} -
         \LC_s\delta\omega_{im}\right)\omega_{jkl}\right]\\
     - \frac{1}{4}\left[\left(\delta\omega_{ij}\delta\omega_{sp} -
         \delta\omega_{sj}\delta\omega_{ip}\right)\omega_{pklm} -
       \left(\delta\omega_{ik}\delta\omega_{sp} -
         \delta\omega_{sk}\delta\omega_{ip}\right)\omega_{plmj}\right]\\
     - \frac{1}{4}\left[\left(\delta\omega_{il}\delta\omega_{sp} -
         \delta\omega_{sl}\delta\omega_{ip}\right)\omega_{pmjk} -
       \left(\delta\omega_{im}\delta\omega_{sp} -
         \delta\omega_{sm}\delta\omega_{ip}\right)\omega_{pjkl }\right]
   \end{multline}
   \qed
 \end{lem}
 Using~\eqref{nat} we get
 \begin{equation}\label{kov}
   \LC_k\delta\omega_{is}\omega_{ism} =
   \Nt_k\delta\omega_{is}\omega_{ism} +
   \frac{1}{6}\delta\omega_{kr}\delta\omega_{rq}\omega_{siq}\omega_{ism}
   = \delta\omega_{kr}\delta\omega_{mr},
 \end{equation}
 since $\Nt\delta\omega \in \Lambda^2_{14}$.  If we multiply~\eqref{ricid1}
 by $\omega_{mkl}$ and use the Bianchi identity as well as~\eqref{kov} we
 obtain~\eqref{ric1}.

 Multiplying~\eqref{ricid2} by $\omega_{mlj}$, and again using the Bianchi
 identity (alternatively: multiply~\eqref{curv} by $\omega_{klm}$), we get
 \begin{equation}
   \label{curv1}
   R_{silr}\omega_{klr} = \left(\LC_s\delta\omega_{ik} -
     \LC_i\delta\omega_{sk}\right) +
   \frac{1}{4}\left(\delta\omega_{ij}\delta\omega_{sp} -
     \delta\omega_{sj}\delta\omega_{ip}\right)\omega_{jpk}.
 \end{equation}
 From~\eqref{curv1} we get that
 \begin{align}
   \label{ro}
   \rho^\star_{km} &= R_{silr}\omega_{klr}\omega_{sim} \\\nonumber &=
   d\delta\omega_{sik}\omega_{sim} - \LC_k\delta\omega_{is}\omega_{sim}
   + \frac{1}{2}\delta\omega_{ij}\delta\omega_{sp}\omega_{jpk}\omega_{sim}.
 \end{align}
 The second term is calculated in~\eqref{kov}.  The last term is manipulated
 using~\eqref{bian} and~\eqref{0}:
 \begin{align}
   \label{cal}
   \delta\omega_{ij}\delta\omega_{sp}\omega_{jpk}\omega_{sim} &=
   \left(-\delta\omega_{pj}\omega_{jki} -
     \delta\omega_{kj}\omega_{jip}\right)\delta\omega_{sp}\omega_{sim}\\
   &= \delta\omega_{sp}\left(\delta\omega_{jp}\omega_{kij}\omega_{sim}
     + \delta\omega_{jk}\omega_{jip}\omega_{sim}\right)\nonumber\\
   \begin{split}&=
     \delta\omega_{sp}\delta\omega_{jp}\left(-\omega_{kjsm}
       + \delta_{ks}\delta_{jm}-\delta_{km}\delta_{js}\right) \nonumber\\
     &\qquad + \delta\omega_{sp}\delta\omega_{jk}\left(-\omega_{jpsm} +
       \delta_{js}\delta_{pm} - \delta_{jm}\delta_{ps}\right)\nonumber
   \end{split}\\
   &= -\norm{\delta\omega}^2\delta_{km} +
   4\delta\omega_{jm}\delta\omega_{jk},\nonumber
 \end{align}
 again, since $\delta\omega \in \Lambda^2_{14}$.  Substituting~\eqref{cal}
 and~\eqref{kov} into~\eqref{ro} we obtain~\eqref{ric2}.
\end{proof}

The equality~\eqref{kov} leads to
\begin{equation}\label{eq}
  d\delta\omega_{sjm}\omega_{sjm}=3\norm{\delta\omega}^2.
\end{equation}
Taking the trace in~\eqref{ric1} and using~\eqref{eq}, we get the formula
for the scalar curvature of a closed $G_2$-structure discovered recently by
Bryant in~\cite{Br1}
\begin{co}
  The scalar curvature of a closed $G_2$-structure is non-positive
  while the $\star$-scalar curvature is non-negative.  These functions
  are given by
  \begin{equation}\label{scal}
    s=-\frac{1}{4}\norm{\delta\omega}^2,\qquad
    s^\star=\frac{1}{2}\norm{\delta\omega}^2.
  \end{equation}
\end{co}
In view of \eqref{scal}, the trace-free part of the Ricci tensor $\rho^0$
has the expression
\begin{equation}\label{staro}
\rho^0=\rho +\frac{1}{28}\norm{\delta\omega}^2g.
\end{equation}
\begin{defn}
  The canonical connection gives us a third Ricci tensor which we denote by
  $\tilde\rho$:
  \begin{equation*}
    \tilde\rho_{ij}=\Rt_{sijs}.
  \end{equation*}
\end{defn}
\begin{co}\label{3co}
  On a $7$-manifold with closed $G_2$-structure the Ricci tensor of the
  canonical connection is related to the Riemannian Ricci tensor through
  the following formula
  \begin{equation*}
    \tilde\rho=\frac{2}{3}\rho.
  \end{equation*}
\end{co}
\begin{proof}
  Taking the trace of~\eqref{curv} we get
  \begin{equation*}
    \rho_{il} = \tilde\rho_{il} - \frac{1}{12}d\delta\omega_{pji}\omega_{pjl}
    + \frac{1}{6}\delta\omega_{is}\delta\omega_{ls}.
  \end{equation*}
  This equality and~\eqref{ric1} completes the proof.
\end{proof}
Furthermore, we have
\begin{pro}\label{form}
  The Ricci three-forms of a closed $G_2$-structure are given by
  \begin{equation}\label{0-1}
    \rho_{jlk}=-\frac{1}{4}\norm{\delta\omega}^2\omega_{jlk}
    - d\delta\omega_{jlk} -
    \frac{1}{2}\left(\delta\omega_{li}\delta\omega_{pi}\omega_{pjk} +
      \delta\omega_{ki}\delta\omega_{pi}\omega_{plj} +
      \delta\omega_{ji}\delta\omega_{pi}\omega_{pkl} \right)
  \end{equation}
  \begin{equation}\label{ro*}
    \rho^\star_{sik} = 2d\delta\omega_{sik}
    +\frac{1}{4}\left(\delta\omega_{ij}\delta\omega_{pj}\omega_{psk} +
      \delta\omega_{kj}\delta\omega_{pj}\omega_{pis} +
      \delta\omega_{sj}\delta\omega_{pj}\omega_{pki}\right).
  \end{equation}
\end{pro}
\begin{proof}
  Substitute~\eqref{lapl} and~\eqref{ro*} into the Weitzenb\"ock
  formula~\eqref{w} to get~\eqref{0-1}.  The cyclic sum in the
  equality~\eqref{curv1} gives~\eqref{ro*}.
\end{proof}

\begin{thm}\label{th1}
  Let $(M^7,\omega)$ be a $G_2$-manifold with closed fundamental form.  If
  $(M,\omega)$ is Einstein and $\star$-Einstein then $M$ is parallel.
\end{thm}
\begin{proof}
  Let $(M,\omega)$ be a $G_2$ manifold with closed $G_2$ structure
  $\omega$. Suppose that both the Einstein and $\star$-Einstein equations
  are satisfied. Proposition~\ref{pro1} in this case yields
  \begin{equation}
    \label{aa}
    \delta\omega_{ij}\delta\omega_{pj} =
    \frac{1}{7}\norm{\delta\omega}^2\delta_{ip}, \qquad
    d\delta\omega_{ijk}\omega_{ijm} =
    \frac{3}{7}\norm{\delta\omega}^2\delta_{km}.
  \end{equation}
  Taking into account~\eqref{aa} and the equalities~\eqref{lapl},
  \eqref{0-1}, we obtain
  \begin{gather}
    \rho^\star_{ijk} = 2d\delta\omega_{ijk} -
    \frac{3}{28}\norm{\delta\omega}^2\omega_{ijk}, \nonumber\\
    \LC^*\LC\omega_{ijk} =
    \frac{1}{7}\norm{\delta\omega}^2\omega_{ijk},\label{la} \\
    \nonumber \rho_{ijk} = -\frac{3}{28}\norm{\delta\omega}^2\omega_{ijk}.
  \end{gather}
  In view of~\eqref{la}, the Weitzenb\"ock formula~\eqref{w} gives
  \begin{equation}\label{fin}
    d\delta\omega = \frac{1}{14}\norm{\delta\omega}^2\omega.
  \end{equation}
  Bryant shows in~\cite{Br1} that on an Einstein manifold with closed
  $G_2$-structure the $\Lambda^3_{27}$-part of $d\delta\omega$ is given by
  the $\Lambda^3_{27}$-part of $\Hodge(\delta\omega\wedge\delta\omega)$.
  Comparing with \eqref{fin}, we see that the $\Lambda^3_{27}$-part of
  $\Hodge(\delta\omega\wedge\delta\omega)$ vanishes.  We need the following
  algebraic
  \begin{lem}\label{l1}
    Let $\alpha$ be a two-form in $\Lambda^2_{14}$.  Then the
    $\Lambda^3_{27}$-part of $\Hodge(\alpha\wedge\alpha)$ vanishes if and
    only if $\alpha=0$.
  \end{lem}
  \begin{proof}
    First note that if $\alpha$ is a two-form in $\Lambda^2_{14}$ then
    $\alpha\otimes \alpha \in S^2\Lambda^2_{14}\subset S^2(\Lambda^2)$.  The
    space of symmetric tensors on $\Lambda^2_{14}$ decomposes as follows
    \begin{equation*}
      S^2\Lambda^2_{14}=V^{77}_{(0,2)}+V^{27}_{(0,2)}+V^1_{(0,0)}.
    \end{equation*}
    Recall that the map $S^2\Lambda^2\to\Lambda^4$ given by $\beta \vee
    \gamma \mapsto \beta \wedge \gamma$ is surjective and equivariant.  By
    Schurs Lemma we may conclude that $\alpha \wedge \alpha \in
    \Lambda^4_{27} \oplus \Lambda^4_1$ if $\alpha \in \Lambda^2_{14}$.

    Now suppose $(\alpha\wedge\alpha)_{\Lambda^4_{27}}=0$.  We may then
    conclude that $\alpha\wedge\alpha=c\Hodge\omega$ for some constant $c$.
    However, as $\alpha$ is a two-form on an odd-dimensional space it is
    degenerate.  Let $X\in\mathbb R^7$ be a non-zero vector such that
    $X\hook\alpha=0$.  Then
    \begin{equation*}
      cX\hook\Hodge\omega = X\hook(\alpha\wedge\alpha) =
      2(X\hook\alpha)\wedge\alpha = 0.
    \end{equation*}
    But the left hand side vanishes only if $c=0$.

    This shows that $\alpha\otimes\alpha\in V^{77}_{(0,2)}$ whenever $\alpha
    \in\Lambda^2_{14}$ satisfies $(\alpha\wedge\alpha)_{\Lambda^4_{27}}=0$.
    But $V^{77}_{(0,2)}$ is precisely the space in which the Ricci-flat
    Riemannian curvature of parallel $G_2$-manifolds takes values, whence
    $\norm{\alpha}^2=0$.
  \end{proof}
  Now, Lemma~\ref{l1} implies $\delta\omega=0$, whence $\LC\omega=0$.
\end{proof}
\begin{rem}
R.L.Bryant observe \cite{Br2} that the following identity holds
\begin{equation}\label{brnew}
\norm{\alpha\wedge \alpha}^2=6\norm{\alpha}^4, \quad \alpha\in \Lambda^2_{14}.
\end{equation}
Clearly \eqref{brnew} implies the lemma.  Note that the constants have
been changed to fit our conventions.
\end{rem}

\section{An integral formula on closed $G_2$ manifold}

Our main technical tool to handle the closed $G_2$-structure on a
compact manifold is the next
\begin{pro}
  Let $(M,\omega,g)$ be a compact $G_2$-manifold with closed
  fundamental form.  Then the following integral formula holds:
  \begin{gather}\label{wnov}
    \int_M\Bigl(\frac{1}{24}\norm{\delta\omega}^4
    +\frac{28}{9}||\rho^0||^2 -\frac{7}{18}
    \rho^0_{pl}\delta\omega_{bl}\delta\omega_{bp}- \frac79(\delta
    R)_{bkl}\omega_{jkl}\delta\omega_{bj}\Bigr)\,dV=0
  \end{gather}
\end{pro}
\begin{proof}
  The first Pontrjagin form $p_1(\nabla)$ of a connection $\nabla$ may
  be defined by
\begin{equation*}
  p_1(\nabla) := \frac{1}{16\pi^2}R_{ijab}R_{klab}e_i \wedge e_j \wedge e_k
  \wedge e_l.
\end{equation*}
The first Pontrjagin class of $TM$ which is the de Rham cohomology
class whose representative element is the first Pontrjagin form of a
some affine connection on $M$ is independent of the connection on $M$.
This implies that $\bigl(p_1(\LC)-p_1(\Nt)\bigr)$ is an exact
$4$-form.  Since the fundamental form $\omega$ is closed, the wedge
product $\bigl(p_1(\LC) - p_1(\Nt)\bigr) \wedge \omega$ is exact.
From Stoke's theorem we obtain
\begin{equation}
  \label{n0}
  \int_M\bigl(p_1(\LC)- p_1(\Nt)\bigr) \wedge \omega = 0.
\end{equation}

However, we may also express the integrand in terms of the curvatures
of $\LC$ and $\Nt$ using that
\begin{equation*}
  16\pi^2p_1(\LC) \wedge \omega = R_{ijab}R_{klab}
  \omega_{ijkl} = R_{abij}R_{klab}\omega_{ijkl},
\end{equation*}
and
\begin{equation*}
  16\pi^2p_1(\Nt) \wedge \omega = \Rt_{ijab}\Rt_{klab}\omega_{ijkl}.
\end{equation*}
From this point on the proof is essentially a brute force calculation
reducing the difference of these two expressions to the
form~\eqref{wnov}.

Since $\Nt$ is a $G_2$-connection we get
\begin{equation}
  \label{n1}
  \Rt_{abij}\omega_{ijkl} = 2\Rt_{abkl}.
\end{equation}
Using~\eqref{curv},~\eqref{0} and~\eqref{n1} we calculate
\begin{multline}
  \label{n2}
  16\pi^2p_1(\LC)\wedge\omega= 2R_{klab}\Rt_{abkl} -
  \frac{2}{3}\bigl(d\delta\omega_{abp} - \LC_p\delta\omega_{ab}\bigr)
  \omega_{pkl}R_{klab}\\
  -\frac{5}{18}\delta\omega_{ai}\delta\omega_{bj}\omega_{ijkl}R_{klab}
  +\frac{8}{9}\delta\omega_{ak}\delta\omega_{bl}R_{klab}.
\end{multline}
Applying~\eqref{curv} to~\eqref{n2} and using~\eqref{0} we obtain after
some calculations that
\begin{multline}
  \label{n3}
  16\pi^2p_1(\LC)\wedge\omega=\\2\Rt_{klab}\tilde
  R_{abkl}+\bigl(d\delta\omega_{abp}-\LC_p\delta\omega_{ab}\bigr)
  \omega_{pkl}\Bigl(-\frac{2}{3}R_{klab}+\frac{1}{3}\tilde
  R_{klab}\Bigr)\\
  +\delta\omega_{ai}\delta\omega_{bj}\omega_{ijkl}\Bigl(
  -\frac{5}{18}R_{klab}+\frac{2}{9}\Rt_{klab}\Bigr)
  +\delta\omega_{ak}\delta\omega_{bl}\Bigl(\frac{8}{9}R_{klab} -
  \frac{1}{9}\tilde R_{klab}\Bigr).
\end{multline}
To calculate the $p_1(\Nt)$ term of the integrand we first observe
that \eqref{curv} implies
\begin{multline}\label{n4}
  \Rt_{ijkl}=\Rt_{klij}+
  \frac{1}{6}\bigl(d\delta\omega_{klp}-\LC_p\delta\omega_{kl}\bigr)
  \omega_{pij}\\
  -\frac{1}{6}\bigl(d\delta\omega_{ijp}-\LC_p\delta\omega_{ij}\bigr)
  \omega_{pkl}
  +\frac{1}{9}\delta\omega_{ks}\delta\omega_{lp}\omega_{spij}-\frac{1}{9}
  \delta\omega_{is}\delta\omega_{jp}\omega_{spkl}.
\end{multline}
Taking into account~\eqref{n1},~\eqref{n4} and~\eqref{0}, we calculate
\begin{multline}
  \label{n5}
  16\pi^2p_1(\Nt)\wedge\omega=2\Rt_{klab}\tilde
  R_{abkl}-\frac{2}{3}\bigl(d\delta\omega_{abp}-\LC_p\delta\omega_{ab}\bigr)
  \omega_{pkl}\Rt_{klab} \\
  - \frac{4}{9}\delta\omega_{ai}\delta\omega_{bj}\omega_{ijkl}\Rt_{klab} +
  \frac{8}{9}\delta\omega_{ak}\delta\omega_{bl}\Rt_{klab}.
\end{multline}
Subtracting~\eqref{n5} from~\eqref{n3} we get
\begin{equation}
  \label{n6}
  16\pi^2\bigl(p_1(\LC)-p_1(\Nt)\bigr)\wedge\omega=A+B+C,
\end{equation}
where $A,B,C$ are given by
\begin{gather}
  \label{n06}
  A=
  \bigl(d\delta\omega_{abp}-\LC_p\delta\omega_{ab}\bigr)
  \omega_{pkl}\Bigl(\frac{1}{3}R_{klab}-\bigl(R_{klab}-
  \Rt_{klab}\bigr)\Bigr)\\ B =
  \delta\omega_{ai}\delta\omega_{bj}\omega_{ijkl}\Bigl(\frac{7}{18}R_{klab} -
  \frac{2}{3}\bigl(R_{klab} - \Rt_{klab}\bigr)\Bigr)\\
  C = \delta\omega_{ak}\delta\omega_{bl}\Bigl(-\frac{1}{9}R_{klab} +
  \bigl(R_{klab}-\Rt_{klab}\bigr)\Bigr).
\end{gather}
We shall calculate each term in~\eqref{n6}.

First observe that \eqref{brnew} implies the useful identity
\begin{equation}\label{brnew1}
 4\delta\omega_{ai}\delta\omega_{bi}\delta\omega_{aj}\delta\omega_{bj}= ||\delta\omega||^4.
\end{equation}
Taking into
account~\eqref{curv},~\eqref{curv1},~\eqref{ric1},~\eqref{kov},
\eqref{cal},~\eqref{eq} and~\eqref{brnew1} we obtain after some
calculation that
\begin{multline}
  \label{n7}
  A= -\frac{19}{24.7}\norm{\delta\omega}^4+
  \frac{1}{3}\norm{d\delta\omega_{abp} - \LC_p\delta\omega_{ab}}^2
  -\frac{8}{3}||\rho^0||^2 +\frac{4}{3}\rho^0_{sp}\delta\omega_{st}\delta\omega_{pt}\\ -
  \frac{1}{9}\bigl(d\delta\omega_{abp} -
  \LC_p\delta\omega_{ab}\bigr)\omega_{ijp}\delta\omega_{ai}\delta\omega_{bj}
  -\frac{1}{9}\bigl(d\delta\omega_{abp} - \LC_p\delta\omega_{ab}\bigr)
  \omega_{pkl}\omega_{smab}\delta\omega_{ks}\delta\omega_{lm}.
\end{multline}
Let $X_s=R_{akbl}\omega_{sak}\delta\omega_{bl}, \quad Y_s=
\delta\omega_{bl}\delta\omega_{li}\delta\omega_{bp}\omega_{ips}$.
Using the first and second Bianchi identity as well as \eqref{curv1},
we get
\begin{align}
  \label{n8}
  R_{klab}\delta\omega_{ak}\delta\omega_{bl} &=
  \frac{1}{2}R_{akbl}\delta\omega_{ak}\delta\omega_{bl}\\
  \nonumber &= -\frac{1}{2} R_{akbl}\LC_s\omega_{sak}\delta\omega_{bl}\\
  \nonumber &= \frac{1}{2}\delta X
  +\frac{1}{2}R_{akbl}\omega_{sak}\LC_s\delta\omega_{bl}\\
  \nonumber &= \frac{1}{2}\delta X + \frac{1}{6}\norm{d\delta\omega}^2 -
  \frac{1}{2}\norm{\LC\delta\omega}^2 +
  \frac{1}{4}\LC_s\delta\omega_{bl}\delta\omega_{li}\delta\omega_{bp}\omega_{ips}\\
  \nonumber
  \begin{split}
    &=\frac{1}{2}\delta X -\frac{1}{4}\delta Y
    +\frac{1}{6}\norm{d\delta\omega}^2-\frac{1}{2}\norm{\LC\delta\omega}^2\\
    &\qquad-\frac{1}{8}||\delta\omega||^4 +\frac{1}{4}
    d\delta\omega_{sil}\delta\omega_{bl}\delta\omega_{bp}\omega_{ips}.
  \end{split}
\end{align}
Applying~\eqref{ric1} to~\eqref{n8}, we obtain
\begin{equation}\label{n9}
  R_{klab}\delta\omega_{ak}\delta\omega_{bl} = \frac{1}{2}\delta X  -
  \frac{1}{4}\delta Y + \frac{1}{6}\norm{d\delta\omega}^2 -
  \frac{1}{2}\norm{\LC\delta\omega}^2 +
  \frac{1}{28}\norm{\delta\omega}^4 -
  \rho^0_{sp}\delta\omega_{st}\delta\omega_{pt}.
\end{equation}
\begin{rem}
  Notice that~\eqref{n9} is the Weitzenb\"ock formula
  \begin{equation*}
    \delta d\delta\omega=\LC^*\LC\delta\omega
    -\frac{1}{14}\norm{\delta\omega}^2\delta\omega +
    \rho^0(\delta\omega, .) + \rho^0(.,\delta\omega)
     +R(\delta\omega)
  \end{equation*}
  for the 2-form $\delta\omega$ on a closed $G_2$-structure.
\end{rem}
Using~\eqref{curv},~\eqref{brnew}
and~\eqref{n9} we get that
\begin{multline}\label{n10}
  C = -\frac{1}{36}\delta(2X - Y) -
  \frac{25}{63.16}\norm{\delta\omega}^4 -
  \frac{1}{54}\norm{d\delta\omega}^2 +
  \frac{1}{18}\norm{\LC\delta\omega}^2 \\ +
  \frac{1}{9}\rho^0_{sp}\delta\omega_{st}\delta\omega_{pt} +
  \frac{1}{6}\bigl(d\delta\omega_{kls} - \LC_s\delta\omega_{kl}\bigr)
  \omega_{sab}\delta\omega_{ak}\delta\omega_{bl}.
\end{multline}

Now we calculate $B$. Denote
$Z_a=R_{klab}\omega_{jkl}\delta\omega_{bj}$.  Using~\eqref{der1}
and~\eqref{curv1} we have the following sequence of equalities
\begin{align}\label{n11}
  \frac{1}{2}R_{klab}\delta\omega_{ai}\omega_{ijkl}\delta\omega_{bj}
  &=R_{klab}\LC_a\omega_{jkl}\delta\omega_{bj}\\&= \nonumber -\delta Z +
  (\delta R)_{bkl}\omega_{jkl}\delta\omega_{bj}-
  R_{klab}\omega_{jkl}\LC_a\delta\omega_{bj}\\\nonumber &= -\delta Z +
  (\delta R)_{bkl}\omega_{jkl}\delta\omega_{bj}
  -\frac{1}{6}\norm{d\delta\omega}^2-
  \frac{1}{2}\norm{\LC\delta\omega}^2 \\\nonumber
  &\qquad+\frac{1}{4}\bigl(d\delta\omega_{abj}-\LC_j\delta\omega_{ab}\bigr)
  \omega_{slj}\delta\omega_{as}\delta\omega_{bl}
\end{align}
Applying~\eqref{curv}, we get
\begin{multline}
  -\frac{2}{3} \bigl( R_{klab}-\tilde
    R_{klab}\bigr)\delta\omega_{ai}\delta\omega_{bj}\omega_{ijkl}=\\
  -\frac{1}{9}\bigl(d\delta\omega_{kls}-\LC_s\delta\omega_{kl}\bigr)
  \omega_{sab}\omega_{ijkl}\delta\omega_{ai}\delta\omega_{bj}\\-
  \frac{2}{27}\delta\omega_{ks}\delta\omega_{lm}\delta\omega_{ai}
  \delta\omega_{bj}\omega_{smab}\omega_{ijkl}.\label{n12}
\end{multline}
Denote $V_k =
\omega_{mab}\delta\omega_{lm}\delta\omega_{bj}\delta\omega_{ai}\omega_{ijkl}$.
Using \eqref{der1}, we obtain
\begin{align}
  \label{n13}
  \frac{1}{2}\delta\omega_{ks}\omega_{smab}\delta\omega_{lm}\delta\omega_{ai}
  \delta\omega_{bj}\omega_{ijkl}=\LC_k\omega_{mab}
  \delta\omega_{lm}\delta\omega_{ai}
  \delta\omega_{bj}\omega_{ijkl}\\= \nonumber
  -\delta V-\frac{3}{2}
  \bigl(d\delta\omega_{klm}-\LC_m\delta\omega_{kl}\bigr)
  \omega_{mab}\omega_{ijkl}\delta\omega_{ai}\delta\omega_{bj}.
\end{align}
Get together~\eqref{n11},~\eqref{n12} and~\eqref{n13} we obtain
\begin{multline}
  B= -\frac1{27}\delta(21Z - 4V) + \frac79(\delta
  R)_{bkl}\omega_{jkl}\delta\omega_{bj}
  -\frac{7}{54}\norm{d\delta\omega}^2- \frac{7}{18}\norm{\LC\delta\omega}^2\\
  +\frac{1}{9} \bigl(d\delta\omega_{klm}-\LC_m\delta\omega_{kl}\bigr)
  \omega_{mab}\omega_{ijkl}\delta\omega_{ai}\delta\omega_{bj}\\
  +\frac{7}{36}\bigl(d\delta\omega_{abj}-\LC_j\delta\omega_{ab}\bigr)
  \omega_{slj}\delta\omega_{as}\delta\omega_{bl}\label{n14}.
\end{multline}
Collecting terms from~\eqref{n7},~\eqref{n10} and~\eqref{n14} we
obtain the
following expression
\begin{multline}
  16\pi^2 \bigl(p_1(\LC) - p_1(\Nt)\bigr)\wedge\omega =
  \delta\Bigl(-\frac{1}{18}X + \frac{1}{36}Y - \frac79Z + \frac{4}{27}V\Bigr)\\
  - \frac{139}{18.56}\norm{\delta\omega}^4 -
  \frac{1}{27}\norm{d\delta\omega}^2 - \frac{8}{3}||\rho^0||^2\\
  + \frac79(\delta R)_{bkl}\omega_{jkl}\delta\omega_{bj} +
  \frac{1}{4}\bigl(d\delta\omega_{kls} -
  \LC_s\delta\omega_{kl}\bigr)
  \omega_{sab}\delta\omega_{ak}\delta\omega_{bl} + \frac{13}{9}
  \rho^0_{pl}\delta\omega_{bl}\delta\omega_{bp}\label{n15}.
\end{multline}
We express the last two terms in a more tractable form.
Applying~\eqref{bian} we have
\begin{align*}
  d\delta\omega_{abp} \delta\omega_{ai} \delta\omega_{bj} \omega_{pij}
  & = - d\delta\omega_{abp} \delta\omega_{ai}\bigl(\delta\omega_{pj}
  \omega_{ibj} + \delta\omega_{ij} \omega_{bpj}\bigr) \\ &= -
  d\delta\omega_{apb} \delta\omega_{ai} \delta\omega_{bj} \omega_{ipj}
  - d\delta\omega_{abp} \delta\omega_{ai} \delta\omega_{ij}
  \omega_{bpj},
\end{align*}
whence
\begin{equation}\label{zaq}
   d\delta\omega_{abp}\delta\omega_{ai}\delta\omega_{bj}\omega_{pij}
   =
   \frac{1}{2}d\delta\omega_{abp}\omega_{bpj}\delta\omega_{ai}\delta\omega_{ji}
   = - 2\rho_{aj}\delta\omega_{ai}\delta\omega_{ji} +
   \frac14||\delta\omega||^4,
\end{equation}
where we used~\eqref{brnew1} and~\eqref{ric1}.

The next equalities are a consequence of \eqref{ric1},~\eqref{brnew}
and~\eqref{kov}
\begin{align}
  \label{zwq}
  \LC_s\delta\omega_{kl}\omega_{sab}\delta\omega_{ak}\delta\omega_{bl}
  &= \delta Y + \frac{1}{4}||\delta\omega||^4 -
  \bigl(d\delta\omega_{sak} - \LC_a\delta\omega_{ks}\bigr)
  \omega_{sab}\delta\omega_{kl}\delta\omega_{bl} \\ \nonumber &=
  \delta Y + 4\rho_{kb}\delta\omega_{kl}\delta\omega_{bl}.
\end{align}
Substitute \eqref{zaq},~\eqref{zwq} into~\eqref{n15},
use \eqref{ric1} and integrate over the compact
$M$  to get
\begin{gather}\label{posl}
  \int_M\Bigl(-\frac{11}{18.28}||\delta\omega||^4 -
  \frac{1}{27}||d\delta\omega||^2 -  \frac{8}{3}||\rho^0||^2 -
  \frac{1}{18}\rho^0_{pl}\delta\omega_{bl}\delta\omega_{bp} + \\
  \nonumber \frac79(\delta R)_{bkl}\omega_{jkl}\delta\omega_{bj}\Bigr)\,dV =
  0.
\end{gather}
We calculate from \eqref{0-1} that
\begin{equation}\label{zuk}
  ||d\delta\omega||^2 = \frac{15}{28}||\delta\omega||^4 +
  12||\rho^0||^2 -
  12\rho^0_{pl}\delta\omega_{bl}\delta\omega_{bp}.
\end{equation}
Substitute~\eqref{zuk} into~\eqref{posl} to obtain~\eqref{wnov}.
\end{proof}
\begin{co}({\bf Integral Weitzenb\"ock formula.})
On a compact $G_2$ manifold with closed fundamental form we have
\begin{equation*}
  \int_MR_{akbl}\delta\omega_{ak}\delta\omega_{bl}\,dV=
  \int_M\Bigl(\frac{1}{3}\norm{d\delta\omega}^2 -
  \norm{\LC\delta\omega}^2 -2\rho_{sp}\delta\omega_{st}\delta\omega_{pt}\Bigr)\,dV.
\end{equation*}
\end{co}
\begin{proof}
The first Bianchi identity implies
$R_{akbl}\delta\omega_{ak}\delta\omega_{bl}=
  2R_{klab}\delta\omega_{ak}\delta\omega_{bl}$.  Apply \eqref{staro} to
\eqref{n9} multiply by two and integrate the obtained equality
over the compact manifold to get the result.
\end{proof}

\section{Proof of the Main Theorem}

We consider the co-differential of the Weyl tensor, $\delta W$ as an
element of $T^*M\otimes\Lambda^2(T^*M)$.  According to the splitting
of the space of two forms, $\delta W$ splits as follows
$$
\delta W=\delta W^2_{14}\oplus\delta W^2_7,
$$
where $\delta W^2_{14}$ is a section of $T^*M\otimes\Lambda^2_{14}(T^*M)\cong
\mathbb R^7\otimes g_2$ while
$\delta W^2_7$ is a section of $T^*M\otimes\Lambda^2_7(T^*M)$.

The Main Theorem is a consequence of the following general
\begin{thm}\label{gen}
Let $(M,\omega,g)$ be a compact $G_2$-manifold with closed
fundamental form $\omega$. Suppose that the $\Lambda^2_7$-part of
the co-differential of the Weyl tensor vanishes, $\delta W=\delta
W^2_{14}$. Then $(M,\omega,g)$ is a Joyce space
\end{thm}
\begin{proof}
  On a 7-dimensional Riemannian manifold the Weyl tensor is expressed
  in terms of the normalized Ricci tensor $h=-\frac{1}{5}\left(\rho
    -\frac{s}{12}g\right)$ as follows
\begin{equation*}\label{w3}
W_{ijkl}=R^g_{ijkl}  -  h_{ik}g_{jl}+ h_{jk}g_{il}-h_{jl}g_{ik}+
h_{il}g_{jk}.
\end{equation*}
The second Bianchi identity implies
\begin{equation}\label{bi2}
\LC_i\rho_{jk}-\LC_j\rho_{ik}=
(\delta R)_{kij}, \quad ds_i=2\LC_j\rho_{ji}, \quad
(\delta W)_{kij}=-4\Bigl(\LC_ih_{jk}-\LC_jh_{ik}\Bigr).
\end{equation}
The condition of the theorem reads
\begin{equation*}
0=(\delta W)_{kij}\omega_{ijl}=\frac{4}{5}\Bigl(\LC_i\rho_{jk}-\LC_j\rho_{ik}\Bigr)\omega_{ijl}
-\frac{2}{15}ds_i\omega_{ikl}.
\end{equation*}
Consequently,
\begin{equation}\label{pr12}
\frac{4}{5}(\delta R)_{kij}\omega_{ijl}\delta\omega_{kl}=
\Bigl(\LC_i\rho_{jk}-\LC_j\rho_{ik}\Bigr)\omega_{ijl}\delta\omega_{kl}=0,
\end{equation}
since $\delta\omega\in\Lambda^2_{14}$.

To apply effectively our integral formula \eqref{wnov} we have to
evaluate one more term.

Denote $K_i=\omega_{ijl}\rho_{jk}\delta\omega_{kl}$. The equation
\eqref{pr12} together with \eqref{ric1} leads to the equality $$
\delta K
=-\frac{1}{2}\rho_{jk}\delta\omega_{jl}\delta\omega_{kl}-2||\rho||^2
$$ which implies, by an integration over the compact $M$, that
\begin{equation}\label{xqaz}
\int_M \rho_{jk}\delta\omega_{jl}\delta\omega_{kl} \,dV = -4\int_M ||\rho||^2\,dV.
\end{equation}
The equalities \eqref{pr12},~\eqref{xqaz},~\eqref{staro}
and~\eqref{wnov} yield
\begin{equation*}
  \int_M
  \Bigl( \frac{1}{24}||\delta\omega||^4 + \frac{14}{3}||\rho^0||^2 \Bigr)\,dV = 0.
\end{equation*}
Hence, Theorem~\ref{gen} follows.
\end{proof}

Clearly, our Main Theorem follows from Theorem~\ref{gen}.\qed

Corollary~\ref{3co} and the main theorem lead to
\begin{thm}
  Any compact $7$-manifold with closed $G_2$-structure which is
  Einstein with respect to the canonical connection is a Joyce space.
\end{thm}

\bibliographystyle{hamsplain}

\providecommand{\bysame}{\leavevmode\hbox to3em{\hrulefill}\thinspace}






\end{document}